\documentclass[11pt]{amsart}
\usepackage{amssymb}
\usepackage{amscd}
\usepackage{color}
\usepackage{url}
\usepackage{hyperref}
\usepackage{xcolor}

\numberwithin{equation}{section}

\DeclareMathOperator*{\Span}{span \;}

\theoremstyle{plain}
\newtheorem{theorem}{Theorem}[section]
\newtheorem{thm}[theorem]{Theorem}
\newtheorem{thm*}{Theorem}

\newtheorem{cor}[theorem]{Corollary}
\newtheorem{lem}[theorem]{Lemma}
\newtheorem{prop}[theorem]{Proposition}
\theoremstyle{definition}

\newtheorem{ex}[theorem]{Example}
\theoremstyle{remark}
\newtheorem{rem}[theorem]{Remark}

\newtheorem*{thCowen}{Theorem (Cowen, 1978)}

\def\CC{\mathbb C}
\def\DD{\mathbb D}
\def\RR{\mathbb R}

\def\AA{\mathcal{A}}
\def\MM{\mathcal{M}}

\def\ds{\displaystyle}
\def\beginpf{\begin{proof}}
\def\endpf{\end{proof}}
\def\beq{\begin{equation}}
\def\eeq{\end{equation}}

\newcommand{\xdownarrow}[1]{%
  {\left\downarrow\vbox to #1{}\right.\kern-\nulldelimiterspace}
}
\def\longdownarrow{\bigg\downarrow}

\begin{document}

\title[Multiplication by a finite Blaschke product]{Multiplication by a finite Blaschke product on weighted Bergman spaces: commutant and reducing subspaces}

\author{Eva A. Gallardo-Guti\'{e}rrez}
\address{Eva A. Gallardo-Guti\'errez \newline
Departamento de An\'alisis Matem\'atico y Matem\'atica Aplicada,\newline
Facultad de Matem\'aticas,
\newline Universidad Complutense de
Madrid, \newline
Plaza de Ciencias 3, 28040 Madrid,  Spain
\newline
and Instituto de Ciencias Matem\'aticas ICMAT,
\newline Madrid,  Spain }
\email{eva.gallardo@mat.ucm.es}

\author{Jonathan R. Partington}
\address{Jonathan R. Partington, \newline
School of Mathematics, \newline
University of Leeds, \newline
Leeds LS2 9JT, United Kingdom}
\email{J.R.Partington@leeds.ac.uk}

\thanks{Both authors are partially supported by Plan Nacional  I+D grant no. PID2019-105979GB-I00, Spain. First author is also supported by
 the Spanish Ministry of Science and Innovation, through the ``Severo Ochoa Programme for Centres of Excellence in R\&D'' (CEX2019-000904-S) and from the Spanish National Research Council, through the ``Ayuda extraordinaria a Centros de Excelencia Severo Ochoa'' (20205CEX001). }


\date{May 2021, revised December 2021 and May 2022}



\begin{abstract}
We provide a characterization of the commutant of analytic Toeplitz operators $T_B$ induced by finite Blaschke products $B$ acting on weighted Bergman spaces which, as a particular instance, yields the case $B(z)=z^n$ on the Bergman space solved recently by by Abkar, Cao and Zhu \cite{ACZ}. Moreover, it extends previous results by Cowen and Wahl in this context and applies to other Banach spaces of analytic functions such as Hardy spaces $H^p$ for $1<p<\infty$. Finally, we apply this approach to study the reducing subspaces of $T_{B}$ in weighted  Bergman spaces.
\end{abstract}


\maketitle

\section{Introduction}
A general question in Operator Theory in Hilbert spaces consists in describing the commutant of a bounded linear operator, since it provides a better  understanding of its action on the whole space. Unless the operator is normal, very little is known about this general question and even for subnormal operators the question is far from being solved completely.

\smallskip

Clearly, such a question has been addressed for particular classes of operators. Indeed, in the pioneering work  \cite{ShW}, Shields and Wallen considered the commutants of a special class of operators
which can be viewed as the coordinate operators on some reproducing kernel Hilbert
spaces of holomorphic functions. Later on, Deddens and Wong \cite{DedWong} studied the communtant of analytic Toeplitz operators $T_f$ on the classical Hardy space $H^2$ whenever the inducing symbol $f$ was univalent. Both works were the starting point of a fruitful line of research where relevant results were proved by Cowen \cite{Cow1, Cow2, Cow3} or Thomson \cite{T1, T2, T3}.

\smallskip

When the underlying space is the Bergman space, Thomson’s theorem states, in particular, that if $f$ is holomorphic on a neighborhood of the closed unit disc $\overline{\mathbb{D}}$, then there is always a finite Blaschke product $B$ and an $H^\infty$-function $g$ on $\mathbb{D}$ such that
$f=g\circ B$ and the commutants of $T_f$ and $T_B$ coincide. Thus, studying the commutant of $T_B$ whenever $B$ is a finite Blaschke provides a clear insight on the commutator of a larger class of functions which turns out to be a difficult question in general.

\smallskip

Likewise, knowing the commutant of $T_B$ provides information about reducing subspaces. Recall that a non-trivial subspace $M$ in $H$ is reducing for a bounded linear operator $T$ if both $M$ and its orthogonal complement $M^\perp$ are invariant under $T$. Equivalently, the orthogonal projection $P_M$ onto $M$ commutes with $T$. At this regard, it is worthy to point out that while the commutant of the analytic Toeplitz operator induced by the monomial $z^n$ for $n\geq 2$ acting on the Hardy space $H^2$ has infinitely many minimal reducing subspaces, such an operator acting on the Bergman space has exactly $n$ distinct minimal reducing subspaces.

\smallskip

Indeed, for Blaschke products of order two, it was shown independently in \cite{SW} and \cite{zhu1} that $T_B$ has exactly
two distinct minimal reducing subspaces in the Bergman space. For finite Blaschke products, there has been an active line of research to determine the conmutant, where the latest progress was due to Douglas, Putinar and Wang \cite{DPW}. We refer to the recent monograph \cite{GH} for relevant results and more on the subject as well as the references therein.

\smallskip

The aim of the present work is twofold. On one hand, we characterize the commutant of analytic Toeplitz operators $T_B$ induced by finite Blaschke products  $B$ on an scale of weighted Bergman spaces, when $B$ is a finite Blaschke product. The case $B(z)=z^n$ was recently solved by Abkar, Cao and Zhu \cite{ACZ}
(following earlier work by \v{C}u\v{c}kovi\'{c} \cite{cuck}),
and we shall later see their result as a special case of our general theorem. Moreover, the characterization extends the results of Cowen and Wahl \cite{CoWa} who proved that on a large scale of Hilbert spaces $\mathcal{H}_k$ of analytic functions on $\mathbb{D}$ satisfying \emph{natural} properties, among them the Bergman space, the commutants of $T_B$ on $\mathcal{H}_k$ and on $H^2$  are the same in the sense that they have the same formulas
in their actions on these spaces, and therefore the operators on these spaces are all
extensions or restrictions of each other.

\smallskip

In addition, as it is well-known, a large class of weighted shift operators $W$ can be realized as the analytic Toeplitz operator of multiplication by $z$ on some weighted Hardy spaces $H^2(\beta)$ of analytic functions on a disk centered at the origin.
It is a classical result due to Shields \cite{shields} that in such a case
every operator in the commutant of $T_z$ on $H^2(\beta)$ is the strong limit of a sequence of polynomials in $T_z$ and thus a member of the multiplier algebra. Accordingly, Abkar, Cao and Zhu's result state that for such weighted shift $W$ operators, the commutant of $W^n$ can be identified with $n$ copies of the multiplier algebra in a natural sense. They  ask to what extent this result can be generalized pointing out that the main difficulty is \cite[Proposition 3]{ACZ}, which is no longer true for general weighted Hardy spaces. From this standpoint, our approach provides such a generalization for analytic Toeplitz operators induced by finite Blaschke products on weighted Bergman spaces.  Indeed, the method works also for other Banach spaces of analytic functions, not necessarily Hilbert spaces, such as the Hardy spaces $H^p$, $1<p<\infty$.

\smallskip

On the other hand, our second aim deals with reducing subspaces for analytic Toeplitz operators induced by finite Blaschke products on weighted Bergman spaces. As we have already pointed out, this is related to study the self-adjoint projections belonging to the commutant $\{T_B\}'$. Nevertheless, by considering appropriate isometries intertwining $T_B$ and the shift operator, it is possible to exhibit reducing subspaces for $T_{B}$.

\smallskip

The rest of this manuscript is organized as follows. In what follows, we introduce the weighted Bergman spaces  where our work takes place, recalling some important issues about them which will be of interest throughout the manuscript. In Section \ref{Section 2} we provide the characterization of the commutant
while in Section \ref{Section 3} we deal with the study of reducing subspaces.

\subsection{The setting} Given a real number $\alpha$, the weighted Bergman space $\AA_{\alpha}$ consists of analytic functions
$f(z)=\sum_{k=0}^{\infty} a_k z^k$ on the unit disc $\DD$ such that its norm
$$
\|f\|_\alpha := \left ( \sum_{k=0}^{\infty}
|a_k|^2 (k+1)^{\alpha} \right )^{1/2}
$$
is finite. Observe that for particular choices of $\alpha$, we recover well-known Hilbert spaces of analytic functions. Namely, for $\alpha = -1$ we have the classical Bergman space $A^2$, $\alpha = 0$ corresponds to the Hardy space $H^2$, and $\alpha = 1$ to the Dirichlet space $\mathcal{D}$. For $\alpha>1$ every function $f\in \AA_\alpha$ admits a continuous extension to the closure of $\DD$.

As introduced, weighted Bergman spaces are particular instances of the so-call weighted Hardy spaces $H^2(\beta)$ where $\beta$ is a sequence of positive numbers. Such spaces were introduced by Shields \cite{shields} in order to study weighted shift operators on Hilbert spaces, though  other classes of operators have been extensively studied on them.

The \emph{multiplier algebra} $\MM_\alpha$ of $\AA_{\alpha}$ consists of analytic functions $\varphi$ on  $\DD$ such that $\varphi\,  f\in \AA_{\alpha}$ for every $f\in \AA_{\alpha}$. It is well known that for $\alpha=-1,0$ the multiplier algebra $\MM_\alpha$ can be identified with $H^\infty$, the space of bounded analytic functions on $\DD$ endowed with the sup-norm. For other instances of $\alpha$ this is no longer the case, and it is an open question to identify them for any $\alpha$.

\smallskip

An important feature which will be needed in what follows is that $\AA_{\alpha}$ spaces admit an expansion in terms of finite dimensional spaces.
More precisely, if $B$ is a finite Blaschke product of degree $n$, namely
\[B(z)= e^{i\theta} \, \prod_{i=1}^{n} \frac{z-\alpha_i}{1-\overline{\alpha_i}z}, \qquad (z\in \mathbb{D})\]
where $\theta\in \RR$ and each $\alpha_i\in \mathbb{D}$ is counted according to its (prescribed) multiplicity, we denote by $K_B$  the
finite-dimensional model space $H^2 \ominus BH^2$. In \cite{CGGP,GPS} the authors give the following theorem:

\begin{thm} \label{thm:CGGPS}
Let $\alpha \in \RR$ and $B$ a finite nonconstant Blaschke product. Then $f \in \AA_{\alpha}$ if and only if $f = \sum_{k=0}^\infty h_k B^k$ (convergence in $\AA_{\alpha}$ norm) with $h_k \in K_B$ and
\begin{equation}\label{expansion norm}
\sum_{k=0}^\infty(k+1)^{\alpha}\|h_k\|_0^2 < \infty.
\end{equation}
\end{thm}

Before proceeding further two remarks are in order. Firstly, in \eqref{expansion norm} we use the
$H^2$ norm $\|\,.\,\|_0$, but since all norms  on finite-dimensional vector spaces are equivalent we could equally use
$\|\,.\,\|_\alpha$.

\smallskip

Secondly, though the result was only given for $\alpha \in [-1,1]$  in \cite{CGGP,GPS}, the argument is based on two facts that
hold for all $\alpha \in \RR$:\\

(1) For $B$ a finite Blaschke product the composition operator $f \mapsto f \circ B$ is bounded on $\AA_{\alpha}$.

\smallskip

\smallskip

(2) For $B$ a finite Blaschke product and  $u \in K_B$ the analytic Toeplitz operator $T_u$ of multiplication by $u$ is bounded on $\AA_{\alpha}$.

\smallskip

These can be deduced from the fact that for any non-negative integer $N$ with $2N+\alpha > -1$, $\AA_{\alpha}$ is equal to the
space of analytic functions $f$ on the unit disc $\DD$ such that
\[
\int_\DD |f^{(N)}(z)|^2 (1-|z|^2)^{2N+\alpha} dA(z) < \infty,
\]
where $dA(z)$ is two-dimensional Lebesgue measure  (see \cite{ZZ,zhu} for more on these spaces).

\smallskip

To prove (1) we use Fa\`a di Bruno's formula
for derivatives of compositions, and for (2) the Leibniz formula for derivatives of products.

\section{A characterization of $\{T_B\}'$ on $\AA_{\alpha}$}\label{Section 2}

In this section, we provide the characterization of $\{T_B\}'$ on $\AA_{\alpha}$ identifying it with $n$ copies of the multiplier algebra in a natural sense, answering the question posed by Abkar, Cao and Zhu in \cite{ACZ}.

Let us consider the space $(\AA_\alpha)^n = \AA_{\alpha} \oplus \cdots \oplus \AA_{\alpha}$ with the $\ell^2$ direct sum norm.
Define the shift $S$ on $(\AA_{\alpha})^n$ by $S(f_1,\ldots,f_n)(z)=(zf_1(z),\ldots,zf_n(z))$. First, we observe the following

\begin{lem}
An operator $A=(A_{jk})_{j,k=1}^n$ on $(\AA_{\alpha})^n$ commutes with $S$ if and only if each $A_{jk}$ lies in the
multiplier algebra of $\AA_{\alpha}$.
\end{lem}

\beginpf
Taking each $f_k$ to be $0$ apart from $k=m$, say, we have
\[
SAf= (SA_{1m}f_m,\ldots, SA_{nm}f_m)^t
\]
and
\[
ASf=(A_{1m}Sf_m,\ldots,A_{nm}Sf_m)^t.
\]
so that each $A_{jm}$ commutes with $S$.

\smallskip

Conversely, it is easy to see by multiplying out matrices, that the condition is
sufficient for $A$ to commute with $S$.
\endpf

Let $u_1,\ldots,u_n$ be an (algebraic)  basis for the model space $K_B$. For example, if $B$ is a finite Blaschke
product with distinct zeros $z_1,\ldots,z_n$, then the functions $u_n(z)=1/(1-\overline z_n z)$ can be used.
For the case $B(z)=z^n$, the functions $u_k(z)=z^{k-1}$, $k=1,\ldots,n$, form a basis of $K_B$.

\smallskip

Now, we have the following commutative diagram:

 \[
 \begin{array}{ccc}
(\AA_{\alpha})^n & {\displaystyle S \atop \displaystyle\longrightarrow}  & (\AA_{\alpha})^n\\
\noalign{\medskip}
 \displaystyle J \big\downarrow & & {\big \downarrow}\displaystyle J \\
\noalign{\medskip}
\AA_{\alpha} &  {\displaystyle\longrightarrow \atop \displaystyle T_B} & \AA_{\alpha}
\end{array}
\]
where the mapping $J$ is defined by
\[
J\left(\sum_{m=0}^\infty a_{1,m} z^m,\ldots, \sum_{m=0}^\infty a_{n,m}z^m\right)
= \sum_{j=1}^n \sum_{m=0}^\infty a_{j,m}  u_j B^m.
\]
Clearly $JS=T_B J$. Moreover $J$ is an isomorphism, since the (squared) norms
\[
\sum_{m=0}^\infty \left\| \sum_{j=1}^n a_{j,m}u_j \right\|^2(m+1)^\alpha \qquad \hbox{and}\qquad
\sum_{m=0}^\infty   \sum_{j=1}^n |a_{j,m}|^2 (m+1)^\alpha
\]
are  equivalent. Accordingly, we have the following

\begin{lem}\label{lem:1st}
Let $\alpha$ be a real number. A bounded linear operator $W$ on $\AA_{\alpha}$ commutes with $T_B$ if and only if it has the
expression $W=JVJ^{-1}$ where
\[
V\left(f_1,\ldots,f_n\right)=\left(\sum_{k=1}^n \phi_{jk} f_k\right)_{k=1}^n
\]
with $(\phi_{jk})_{j,k=1}^n \subset \MM_\alpha$.
\end{lem}

 From Theorem \ref{thm:CGGPS}, we have a unique decomposition
$f= \sum_{j=1}^n u_j f_j(B)$ for $f \in \AA_{\alpha}$, where $f_1,\ldots,f_n$ also lie in $\AA_{\alpha}$.
Using this, we can obtain an explicit formula from Lemma \ref{lem:1st}.

\begin{thm}\label{thm:w}
Let $\alpha\in \RR$.
A bounded linear operator $W$ on $\AA_{\alpha}$ commutes with $T_B$ if and only if
\beq\label{eq:w}
W\left(\sum_{j=1}^n u_j f_j(B) \right)=\sum_{k=1}^n \phi_k f_k(B)
\eeq
for some $\phi_1,\ldots,\phi_n$ in $\MM_\alpha$.
\end{thm}

\beginpf
We see first that every expression as in \eqref{eq:w} defines an operator that
commutes with $T_B$.\\

For the converse,  we have $W=JVJ^{-1}$ for some
$V$ commuting with $S$, as above, and
\begin{eqnarray*}
JVJ^{-1}\left(\sum_{j=1}^n u_j f_j(B) \right) &=& JV (f_1,\dots,f_n) \\
&=& J \left(\sum_{k=1}^n \phi_{jk} f_k\right) \\
&=& \sum_{k=1}^n \sum_{j=1}^n u_j\phi_{jk}(B) f_k(B).
\end{eqnarray*}
Writing $\phi_k=\sum_{j=1}^n u_j \phi_{jk}(B)$, we see easily that
$\phi_k \in \MM_\alpha$, since each $u_j$ is a rational function and also a multiplier of $\AA_{\alpha}$.
\endpf

From this we have the special case $B(z)=z^n$ given in \cite{ACZ}.

\begin{cor}
Let $\alpha\in\RR$ and let $n$ be a positive integer. A bounded linear operator $W$ on $\AA_{\alpha}$ commutes with
$T_{z^n}$ if and only if, on
decomposing a function $f \in \AA_{\alpha}$ uniquely as
\[
f(z)= \sum_{j=0}^{n-1} z^j f_j(z^n),
\]
with $f_0,f_1,\ldots, f_{n-1}$ in $\AA_{\alpha}$, we have
\beq\label{eq:wzn}
W\left(\sum_{j=0}^{n-1} z^j f_j(z^n) \right)=\sum_{k=0}^{n-1} \phi_k f_k(z^n)
\eeq
for some $\phi_0,\ldots,\phi_{n-1}$ in $\MM_\alpha$.
\end{cor}
\beginpf
This follows from Theorem \ref{thm:w} since $\{1,z,z^2\ldots,z^{n-1}\}$ is a basis for $K_{z^n}$.
\endpf

\subsection{A remark on a theorem of Cowen and Wahl}

As we noted in the Introduction, Cowen and Wahl \cite{CoWa} by means of a generalization of a theorem of Cowen \cite{Cow1} proved that on a large scale of Hilbert spaces $\mathcal{H}_k$ of analytic functions on $\mathbb{D}$, among them the Bergman space, the commutants of $T_B$ on $\mathcal{H}_k$ and on $H^2$  are the same in the sense that they have the same formulas in their actions on these spaces.

In particular, for $\alpha \in \DD$ and $k_\alpha$ denoting the reproducing kernel for $H^2$, Cowen showed the following result.

\begin{thCowen}
Let $T$ be a linear bounded operator acting on the Hardy space $H^2$ and $f\in  H^\infty$. The following are equivalent:
\begin{enumerate}
\item $T$ commutes with $T_f$.
\item For all $\alpha \in \mathbb{D}$, $T^\ast k_{\alpha} \perp (f-f(\alpha))H^2$.
\item There exists $J\subset \mathbb{D}$ such that $\sum_{\alpha \in J} (1-|\alpha|)=\infty$ and for all $\alpha\in J$,  $T^\ast k_{\alpha} \perp (f-f(\alpha))H^2$.
\end{enumerate}
\end{thCowen}

In the case where $f=B$, a finite Blaschke product, we may show easily that Cowen's conditions are
equivalent to our condition \eqref{eq:w}.

In particular, by taking $W=T$ we have that for every $g \in H^2$,
\begin{eqnarray*}
\langle W^* k_\alpha, (B-B(\alpha)g) \rangle &=& \overline{W((B-B(\alpha))g)(\alpha)}
\end{eqnarray*}
Now, assuming \eqref{eq:w}, we obtain
that $W(Bg)=BW(g)$ so that
\[
W(Bg)(\alpha)=B(\alpha)W(g)(\alpha)=W(B(\alpha)g)(\alpha)
\]
for all $\alpha \in \DD$,
 from which
we see that $\langle W^* k_\alpha, (B-B(\alpha)g) \rangle=0$.\\

Conversely, suppose that
$$W(Bg-B(\alpha)g)(\alpha)=0$$
for all $\alpha \in J \subset \DD$ where $J$ is a non-Blaschke subset of $\DD$,
and consider $W(u_j f_j(B))$ where we use the notation
of Theorem \ref{thm:w}, so that $u_j \in K_B$ is fixed and $f_j \in H^2$ is arbitrary.

We begin by identifying $W(u_j B^n)$ for $n=0,1,2,\ldots$. Observe that
\[
W(u_j B^{n+1})(\alpha)= B(\alpha) W(u_j B^n)(\alpha)
\]
 for all $\alpha \in J$.

But if two functions $h_1,h_2 \in H^2$ satisfy $h_1(\alpha)=B(\alpha) h_2(\alpha)$ for all $\alpha \in J$, then
we clearly have $h_1=Bh_2$; thus in our case,
$W(u_j B^{n})= B^{n}W(u_j)$ for all $n=0,1,2,\ldots$ and we
conclude by continuity and linearity that
$W(u_j f(B))=\phi_j f(B)$ for some $\phi_j \in H^\infty$. Accordingly, \eqref{eq:w} follows.

\medskip

Analogously, the Cowen--Wahl Theorem can be derived as an application of Theorem \ref{thm:w}.

\subsection{On $\{T_B\}'$ on $H^p$ spaces}

The previous approach may be considered as well on other spaces of analytic functions, not necessarily Hilbert spaces. This is the case of
the classical Hardy spaces $H^p$ for $1 < p < \infty$. Recall that an analytic function $f$ on $\mathbb{D}$ belongs to $H^p$ if
the norm
$$
\|f\|_p=\left ( \sup_{0\leq r<1} \int_{0}^{2\pi}
|f(re^{i\theta})|^p \, \frac{d\theta}{2\pi}\right )^{1/p}
$$
is finite.  A classical result due to Fatou (see \cite{Du}, for instance) states that the radial limit $f(e^{it}):= \lim_{r\to 1^{-}} f(re^{it})$ exists a.e. on the boundary $\mathbb{T}$ of the unit disc.
In this regard, it is well known that $H^p$ can be regarded as a closed subspace of $L^p(\mathbb{T})$. In particular, for $p=2$ the Lebesgue space  $L^2(\mathbb{T})$ may be decomposed as the orthogonal sum
\begin{equation}\label{eq0}
L^2(\mathbb{T})=H^2 \oplus \overline{H^2_0},
\end{equation}
where $\overline{H^2_0}=\{ f\in L^2(\mathbb{T}): \overline{f}\in H^2  \mbox{ and } f(0)=0\}$.  Note that in the identity \eqref{eq0} we
are identifying $H^2$ through the non-tangential boundary values of the $H^2$ functions.

Let $B$ be a finite Blaschke product of degree $n$ and denote by $K_B$ the associated \emph{model space} in $H^p$, namely,
$$H^p \cap B\overline{H^p_0},$$
which turns out to be an $n$-dimensional space which does not depend on $p$.

For $p=2$, the Wold Decomposition Theorem yields that
$$H^2= K_B \oplus B K_B \oplus B^2 K_B \oplus \ldots,$$
as an orthogonal decomposition. Nevertheless, in the case of general $p$, we take an orthonormal basis $u_1,\ldots,u_n$ of
$K_B$ (in the $H^2$ sense) and we claim that
\[
f \in H^p \iff f=\sum_{k=1}^n u_k f_k(B) \qquad \hbox{for some} \quad f_1,\ldots,f_n \in H^p.
\]
Moreover, there is an equivalence of norms between $\|f\|_p$ and $\sum_{k=1}^n \|f_k\|_p$.

Since the composition operator $C_B$ is  bounded we have,
for finite sums,
\[
\left\| \sum_{j=0}^m \sum_{k=1}^n a_{jk} u_k B^j \right\|_p \le \sum_{k=1}^n \|u_k\|_\infty
\left\| \sum_{j=0}^m a_{jk} B^j  \right\|_p
\le \max \|u_k\|_\infty\sum_{k=1}^n  \left\|\sum_{j=0}^m a_{jk} z^j \right\|_p.
\]
This means that there is a continuous  mapping $H^p \times \ldots \times H^p \to H^p$ given by
\begin{equation}\label{eq:embedding}
(f_1,\ldots,f_n) \mapsto \sum_{k=1}^n u_k f_k(B),
\end{equation} with dense range (since the union of the subspaces $B^jK_B$, $j=0,1,2,\ldots$  is dense).

Using the isomorphic duality $(H^p)^*\simeq H^q$ induced  by the $H^2$ inner product
we see that in fact the mapping \eqref{eq:embedding} is also bounded below, and so it has closed range and is a surjection.
Accordingly, the identification of the commutant of $T_B$ on $H^p$ provided by the formula in Theorem \ref{thm:w} also holds.

\begin{rem}
Observe that in the previous discussion it is not necessary that $u_1,\ldots,u_n$ to be orthonormal, simply an algebraic basis, since in such a case, the duality argument would involve taking the biorthogonal basis $v_1,\ldots,v_n$ and the duality formula
\[
\left\langle \sum_{k=1}^n u_k f_k(B) , \sum_{k=1}^n v_k g_k(B) \right \rangle_{H^p \times H^q}= \sum_{k=1}^n \langle f_k,g_k \rangle_{H^p \times H^q}.
\]
\end{rem}

\section{Minimal reducing subspaces of $T_B$  on the Bergman space}\label{Section 3}

As it was pointed out in the Introduction, recently there has been an increasing interest on reducing subspaces of multiplication operators whenever they act on spaces of analytic functions, and in particular on the Bergman space $A^2$. Clearly the works by Douglas and coauthors \cite{DSZ}, \cite{DPW} or those by Guo and Huang \cite{GH1}, \cite{GH2}, \cite{GH3} are good examples in this regard.

\smallskip

In this setting, if $M$ is a reducing subspace of the multiplication operator $T_B$, the orthogonal projection $P_M$ onto $M$ commutes with both $T_B$ and its adjoint $T_B^\ast$. Accordingly, if $\mathcal{W}(B)$ denotes the von Neumann algebra generated by the multiplication $T_B$ on $A^2$, that is, the unital $C^*$-algebra on $A^2$ which is closed in the weak operator topology, it is plain that $P_M$ belongs to the commutant $\mathcal{W}(B)'$. Note that, in particular, $\mathcal{W}(B)'=\{T_B, T_B^\star\}'$.

\smallskip

As   noted in \cite{GH}, on the Hardy space every non-trivial closed invariant subspace $N$ in $H^2\ominus BH^2$ satisfies that $B^m N \perp B^n N$ for $n\neq m$, and hence the direct sum of all $B^k N$, $k\geq 0$  gives a reducing subspace of $T_B$. On the
Bergman space, in \cite{HSXY}, the authors showed that there is always a canonical reducing subspace $M_0$ such that $T_B$ restricted on $M_0$ is unitarily equivalent to the Bergman shift. Their main tool relies on the Hardy space of the bidisk.
We shall provide some similar results in the context of weighted Bergman spaces.

\smallskip

On the other hand, recently Douglas, Putinar and Wang \cite{DPW} proved that the von Neumann algebra $\mathcal{W}(B)'$ is abelian for any finite Blaschke product $B$ whenever $T_B$ acts on the Bergman space. Indeed, a clear consequence in this context is the following

\begin{prop}
For every finite Blaschke product, $T_B$ has a non-trivial closed hyperinvariant subspace in the Bergman space.
\end{prop}

Note that if $A \in \{T_B\}'$, having into account that $T_B$ has a non-trivial reducing subspace $M$, the orthogonal projection onto $M$, $P_M$, belongs to $\{T_B\}'$ and hence, $AP_M = P_M A$.  So the subspace $M$ is invariant for every operator commuting with $T_B$.

\medskip

In order to find minimal reducing subspaces of $T_B$ it is helpful to
understand the minimal rank-1 projections that commute with $T_B$.

Since $1 \in \AA_{\alpha}$ we clearly have that every multiplier $\phi \in M_\alpha$ corresponds to a function
in $\AA_{\alpha}$. We therefore carry on from Theorem \ref{thm:w} and write $\phi_k  = \sum_{j=1}^n u_j \phi_{jk}(B)$,
where the $\phi_{jk}$ lie in $\AA_{\alpha}$.

Then \eqref{eq:w} can be written as
\[
W\left(\sum_{j=1}^n u_j f_j(B) \right)=\sum_{k=1}^n \sum_{j=1}^n u_j \phi_{jk}(B)  f_k(B),
\]
or, in an obvious notation
\[
\begin{pmatrix}f_1 \\ f_2 \\ \cdots \\  f_n \end{pmatrix}
\mapsto
\begin{pmatrix}
\phi_{11} & \phi_{12} & \ldots & \phi_{1n} \\
\phi_{21} & \phi_{22} & \ldots & \phi_{2n} \\
\ldots & \ldots & \ldots & \ldots \\
\phi_{n1} & \phi_{n2} & \ldots &  \phi_{nn} \end{pmatrix}
 \begin{pmatrix}f_1 \\ f_2 \\ \cdots \\  f_n \end{pmatrix} = \Phi \begin{pmatrix}f_1 \\ f_2 \\ \cdots \\  f_n \end{pmatrix}, \quad \hbox{say} .
\]
Thus the space of  operators commuting with $T_B$ is isomorphic to the space of $n \times n$ matrices of
multipliers of $\AA_{\alpha}$. (Of course, to get a natural isometry one would need to work with an equivalent
norm on $\AA_{\alpha}$.)

However, algebraically at least, a minimal projection $P$ such that $P^2=P$, corresponds to a rank-one matrix
$\Phi$ such that $\Phi^2=\Phi$. The following example  exhibits some rank-one projections in the commutant of $T_{z^2}$.

\begin{ex}
Let $B(z)=z^2$, so that every $f \in \AA_{\alpha}$ decomposes as
$f(z)=f_1(z^2)+ zf_2(z^2)$. This is mapped to the function
\[
\phi_{11}(z)f_1(z^2) +\phi_{12}(z)f_2(z^2) + z \phi_{21}(z) f_1(z^2) + z \phi_{22}(z) f_2(z^2),
\]
where the $\phi_{ij}$ are in $M_\alpha$.

Some  examples of rank-1 idempotents are $\ds \begin{pmatrix}1 & 0 \\ 0 & 0\end{pmatrix}$
and $\ds \begin{pmatrix}1 & 0 \\ p(z) & 0\end{pmatrix}$, where $p$ is a polynomial, or indeed any matrix with rank 1 and trace equal to 1,
such as $\ds \begin{pmatrix}\frac12  & \frac12 \\ \frac12 & \frac12\end{pmatrix}$.
\end{ex}

In the Hardy space the decomposition in Theorem\ref{thm:CGGPS} is orthogonal, a fact which we have exploited already. This
is not the case in the Bergman space, but next we propose an alternative.

\medskip

\subsection{An orthogonal decomposition in the Bergman space $\mathcal{A}_{-1}$} Note that if $B$ is a finite Blaschke product of degree $N$, then the multiplication operator $T_B$ is not an isometry on $\mathcal{A}_{-1}$ and in the decomposition of Theorem \ref{thm:CGGPS} we do not have $K_B$   orthogonal to $B\mathcal{A}_{-1}$ in general.
However, $T_B$ is an injective operator whose range is closed with finite codimension (it is Fredholm),
and moreover it is bounded below.
Thus we can construct an orthogonal decomposition
\beq\label{eq:Gdecomp}
\mathcal{A}_{-1} = K_0 \oplus BK_1 \oplus B^2K_2 \oplus \ldots
\eeq
as follows. We take $K_0$ to be the orthogonal complement of the closed subspace $B \mathcal{A}_{-1}$, a finite-dimensional
subspace which, in the case of distinct zeros, is spanned by Bergman reproducing kernels. We then
take $K_0+BK_1$ to be the orthogonal complement of $B^2\mathcal{A}_{-1}$ (the span of reproducing kernels and derivative kernels)
 in such a way that $BK_1$ is orthogonal to
$K_0$. We continue in this way, and for each $k$ the space $K_k$ has dimension equal to $N$, the
degree of $B$.

\smallskip

If $B$ is a monomial, $B(z)=z^N$, such a decomposition provides a straightforward description of  reducing subspaces as follows.
For each $0\leq j\leq N-1$, let
$$
\mathcal{M}_j=\overline{\Span \{z^{j+kN}:\; k\geq 0 \}}^{\mathcal{A}_{-1}}.
$$
Clearly, each $\mathcal{M}_j$ is a reducing subspace for $T_{z^N}$ and Stessin and Zhu proved that they are the only minimal reducing subspaces for $T_{z^N}$ in the Bergman space. Moreover, $T_{z^N}$ has exactly $2^N$ distinct reducing subspaces counting the trivial ones (see \cite[Theorem B]{SZ}).

\medskip

If $\displaystyle B(z)= \left( \frac{z-a}{1-\overline a z}\right)^N$ for $a \in \DD, a \ne 0$, taking
$$K_0=\Span \{z^j/(1-\overline a z)^{j+2}:\;  0 \le j \le N-1\}$$
namely, to within constants, the successive derivative kernels, we have $K_0$  is orthogonal to $B\mathcal{A}_{-1}$.
Next, $BK_0$ is orthogonal to both $K_0$ and $B^2 \mathcal{A}_{-1}$, since the span of $K_0$ and $BK_0$
is the same as
$$\Span \{z^j/(1-\overline a z)^{j+2}: 0 \le j \le 2N-1\}.$$
Arguing similarly, we deduce that $B^i K_0$ is also orthogonal to $B^{i-1}K_0$ and $B^i \mathcal{A}_{-1}$ for every $i\geq 1$. Accordingly,
\eqref{eq:Gdecomp} turns out to be
\beq
\mathcal{A}_{-1} = K_0 \oplus BK_0 \oplus B^2K_0 \oplus \ldots
\eeq
and the following statement follows as a consequence:

\begin{prop}
Let $\displaystyle B(z)= \left( \frac{z-a}{1-\overline a z}\right)^N$ for $a \in \DD, a \ne 0$. For each $0\leq j\leq N-1$, the subspaces
$$
\mathcal{M}_j=\overline{\Span \left \{ \Big (\frac{z}{(1-\overline a z)^2}\Big )^{j+kN}:\; k\geq 0 \right\}}^{\mathcal{A}_{-1}}.
$$
are reducing for $T_B$ in $\mathcal{A}_{-1}$.
\end{prop}

The main theorem in \cite{HSXY} (see, more recently, \cite{ghosh}
for a multi-dimensional variant) asserts that if
$B$ is a finite Blaschke product of degree $N > 1$, then the analytic Toeplitz operator $T_B$ on the Bergman space is   reducible and there is a reducing subspace on which the restriction of $T_B$   is unitarily equivalent to the
shift $S=T_z$  on $\AA_{-1}$. In the scale of weighted Bergman spaces we have the following.

\begin{theorem}\label{unitarily equivalent} For $\alpha \in \RR$ and $B(z)=z^n$ for an  integer $n \ge 1$ the space $\AA_\alpha$ has a closed $T_B$-reducing subspace on which $T_B$ is unitarily equivalent to the
shift $S=T_z$ on $\AA_\alpha$. For a general finite Blaschke product $B$  there is an
equivalent Hilbert space norm $\|\cdot\|_B$ on $\AA_\alpha$ under which $T_B$ has a reducing subspace
on which it is unitarily equivalent
to the shift on $\AA_\alpha$.
\end{theorem}
\beginpf
For $B(z)=z^n$ we have an obvious reducing subspace $V_B$
spanned by $z^{n-1},z^{2n-1},\ldots$ (its orthogonal complement is clearly invariant under $T_B$),
and since $\|z^{mn-1}\|^2_\alpha = (mn)^{\alpha}$ we see that we have the following commutative diagram,
where $J$ is the unitary mapping with $J(z^k)=z^{(k+1)n-1}/n^{\alpha/2}$ for $k=0,1,2,\ldots$.
\[
\begin{array}{ccc}
&S& \\
\AA_\alpha & \longrightarrow & \AA_\alpha \\
J\,\longdownarrow && \longdownarrow\, J \\
V_B & \longrightarrow & V_B\\
&T_B&
\end{array}
\]
For a general $B$ we return to the decomposition in Theorem \ref{thm:CGGPS}.
Writing $f = \sum_{k=0}^\infty h_k B^k$ (convergence in $\AA_{\alpha}$ norm) with $h_k \in K_B$ we
define
\[
\|f\|_B^2=\sum_{k=0}^\infty(k+1)^{\alpha}\|h_k\|_0^2 < \infty.
\]
We may choose $h \in K_B$ to have $\AA_\alpha$ norm 1, and then let
the subspace $V_B$ be spanned by $h, hB, hB^2, \ldots$, noting that $\|h B^k\|^2_B =  (k+1)^\alpha$,
and hence we have the same commutative diagram as before with $J(z^k)=h B^k$.
\endpf

Observe that it we endow $V_B$  with the usual $\AA_\alpha$ norm instead, then the operator $J$ is merely an isomorphism,
and in general not unitary.

\begin{rem}
In \cite{GPS} the authors showed that in a scale of weighted Dirichlet spaces including the Bergman space $\AA_{-1}$, for
every finite Blaschke product $B$, it is always possible to renorm the space (with an
equivalent norm) such that $B$ enjoys the \emph{wandering subspace property} (see \cite[Theorem 3.3]{GPS}). Accordingly, both Theorem 3.3 in \cite{GPS} and Theorem \ref{unitarily equivalent} show, in particular, that the geometry of the space plays a significant role in order to study such questions since the answer depends on the norm expression.
\end{rem}

\subsection{A final remark: An orthogonal decomposition in the Bergman space related to the $T_B$ action}
\mbox{ }

As a final observation, we note that most of the previous results relied on \emph{ad-hoc} orthogonal decompositions
of the underlying space which allowed us to simplify strategies. At this regard, we finish by exhibiting another orthogonal decomposition of the Bergman space $\mathcal{A}_{-1}$ which is closer to the action of $T_B$ as a shift.

\smallskip

Let $B$ be a Blaschke product of degree $N$. Observe that $T_B^k \mathcal{A}_{-1}$ is a closed subspace of $\mathcal{A}_{-1}$ of codimension $Nk$ for $k=1,2,\ldots$.

\smallskip

Let us write $X_k= T_B^k  \ominus T_B^{k+1} \mathcal{A}_{-1}$ for $k=0,1,2, \ldots$, where $T_B^0$ denotes the identity operator. Note that each  $X_k$ has dimension $N$. We therefore have
\begin{equation}\label{eq:xdecomp}
\mathcal{A}_{-1} = X_0 \oplus X_1 \oplus X_2 \oplus \ldots,
\end{equation}
which is an orthogonal decomposition. Clearly, it allows one to exhibit suitable orthonormal bases of $\mathcal{A}_{-1}$ by combining orthonormal bases of the $X_k$.

Observe that, in particular,
\[
T_B^k \mathcal{A}_{-1} = X_k \oplus X_{k+1} \oplus \ldots 
\]
for $k=0,1,2, \ldots$ Accordingly, $T_B$ acts as a shift respect to the orthogonal decomposition \eqref{eq:xdecomp}.

On the other hand, note that if $W$ commutes with $T_B$ then for a vector $x=B^ky$ in $T_B^k \mathcal{A}_{-1}$ we have
$$Wx = WB^k y = B^k Wy \in T_B^k  \mathcal{A}_{-1}.$$
In particular, this implies that $W$ has a lower triangular block matrix with respect to the decomposition
\eqref{eq:xdecomp}. Moreover, if $W$ is also self-adjoint then the matrix is block-diagonal with
every block self-adjoint.

Consequently, every projection in $\{T_B\}'$ admits a matrix representation respect to \eqref{eq:xdecomp}
which is block-diagonal with every block self-adjoint. Such matrix representation may be of help in understanding
the reducing subspaces for $T_B$.


\end{document}